\documentclass[11pt,a4paper]{amsart}
\usepackage{amsfonts,euscript,verbatim}

\theoremstyle{plain}
\newtheorem{thm}{Theorem}
\newtheorem{prop}[thm]{Proposition}

\newtheorem{lem}[thm]{Lemma}
\theoremstyle{definition}

\newtheorem{defi}[thm]{Definition}

\def\limsup{\operatornamewithlimits{\overline{\text{\rm{lim}}}}}
\def\liminf{\operatornamewithlimits{\underline{\text{\rm{lim}}}}}
\def\cvto  {\operatornamewithlimits{\to}}

\def\A{\EuScript{A}}
\def\B{\EuScript{B}}
\def\F{\EuScript{F}}
\def\M{\EuScript{M}}

\def\RR{\mathbb{R}}

\def\NN{\mathbb{N}}

\def\Z{\mathcal{Z}}

\def\glu{\!\!\!\!\!\!\!\!}
\newcommand{\eps}{\varepsilon}
\newcommand{\e}{{\mathrm{e}}}
\newcommand{\eqdef}{\stackrel{\scriptscriptstyle\rm def}{=}}
\DeclareMathSymbol{\varnothing}{\mathord}{AMSb}{"3F}
\DeclareMathOperator{\diam}{diam}

\DeclareMathOperator{\osc}{osc}
\DeclareMathOperator{\var}{var}
\renewcommand{\emptyset}{\varnothing}

\begin{document}
\date{\today}
\title{recurrence spectrum in smooth dynamical systems.}
\begin{abstract}
We prove that for conformal expanding maps the return time does have constant
multifractal spectrum.
This is the counterpart of the result by Feng and Wu in the symbolic setting.
\end{abstract}

\author{Beno\^\i t Saussol}
\address{CNRS UMR 6140 - LAMFA, Universit\'e de Picardie
Jules Verne, 33 rue Saint Leu, 80039 Amiens, France}
\email{benoit.saussol@u-picardie.fr}
\urladdr{http://www.mathinfo.u-picardie.fr/saussol/}

\author{Jun Wu}
\address{Department of Mathematics and Nonlinear Science Center, Wuhan University,
Wuhan, 430072, People's Republic of China}
\email{wujunyu@public.wh.hb.cn}


\maketitle

\section{Introduction}
Let $f:\M\to \M$ be a map on a smooth compact manifold $\M$ with metric $d$.
Let $\Lambda\subset \M$ be a compact $f$-invariant set.
For $x\in\M$ we can define the return times in $r$-neighborhoods by
\[
\tau_r(x) = \inf\{n>0\colon d(f^nx,x)<r\}.
\]
We put $\tau_r(x)=\infty$ if the orbit of $x$ does not come back close to $x$.

When $\M$ is an interval, $f$ is a piecewise monotonic transformation
with some regularity (piecewise $C^{1+\alpha}$ suffices) and $\mu$ is an
$f$-invariant ergodic measure with nonzero entropy,
Saussol, Troubetzkoy  and Vaienti \cite{stv}, building on results
by Ornstein and Weiss \cite{ow} and Hofbauer and Raith \cite{h,hr}, proved that
\begin{equation}\label{eq:stv}
\lim_{r\to 0} \frac{\log\tau_r(x)}{-\log r} = \dim_H\mu
\quad\text{for $\mu$-almost every $x$}
\end{equation}
where $\dim_H\mu$ stands for the Hausdorff dimension of the measure $\mu$.
Barreira and Saussol established, with completely different techniques,
Equation~\eqref{eq:stv} for equilibrium states
preserved by repellers \cite{bs2} and Axiom~A diffeomorphisms \cite{bs1} in any dimension.

{}From now on we assume that $\Lambda\subset \M$ is a \emph{conformal repeller} of a $C^{1+\alpha}$
map $f\colon \M\to\M$. 
That is to say that $\Lambda$ is a compact forward invariant set
and that $f$ is $C^{1+\alpha}$ on $\M$, expanding and conformal on $\Lambda$. 
Moreover we suppose that there exists an open set $V\subset \M$ such that 
$\Lambda=\bigcap_{n\ge 0} f^{-n} V$. 

{}From the multifractal analysis of conformal expanding maps it is well known that
for any real $s\in [0,\dim_H\Lambda]$ there exists an ergodic equilibrium state $\mu_s$
(in fact there are many such measures for non extremal values of $s$).
Therefore Equation~\eqref{eq:stv} together with basic properties of dimension of measures
imply that
\begin{equation}\label{eq2}
\dim_H \left\{x\in \Lambda\colon \lim_{r\to0}\frac{\log\tau_r(x)}{-\log r}
\text{ exists and equals } s \right\} \ge \dim_H\mu_s = s.
\end{equation}
The reader familiar with multifractal analysis will believe that in fact the
dimension of this set is larger than this, and at this point may expect
a nice non-trivial (strictly) concave analytic curve.
In fact, this is really not the case, as our first theorem says

\begin{thm}\label{thm:1}
Let $\Lambda\subset\M$ be a conformal repeller of the $C^{1+\alpha}$ map $f\colon \M\to \M $.
For any $\alpha\le\beta$ in $[0,\infty]$ we have
\[
\dim_H \left\{x\in \Lambda\colon
\liminf_{r\to0}\frac{\log\tau_r(x)}{-\log r} = \alpha
\text{ and }
\limsup_{r\to0}\frac{\log\tau_r(x)}{-\log r} = \beta
\right\} = \dim_H \Lambda.
\]
\end{thm}
Theorem~\ref{thm:1} implies in particular that for all values of $s\in [0,\infty]$
the dimension in the left hand side of Equation~\eqref{eq2} equals $\dim_H \Lambda$.

The corresponding result of Theorem~\ref{thm:1} in the symbolic case was established
by Feng and Wu \cite{fw}. Although our method relies on some essential ideas present
in the above mentioned paper, we emphasize that our proof is not a direct translation
of the symbolic result.
\footnote{After the completion of this work we learn that L.~Olsen has
obtained --- independently and with different techniques --- a similar
result for conformal Iterated Function Systems satisfying some strong
separation condition, which is equivalent to a conjugacy to a subshift
of finite type, the case considered in Theorem~\ref{thm:2}.}

In Section~\ref{sec:osc} we prove the theorem in the situation where
the system $(f,\Lambda)$ is conjugated to a subshift of finite type.
It happens that our problem can be reduced to this case and this is proven
in Section~\ref{sec:no-osc}.

\tableofcontents

\section{System conjugated to a subshift of finite type}~\label{sec:osc}

In this section we give a proof of Theorem~\ref{thm:1} in a special case.

\begin{thm}\label{thm:2}
Let $\Lambda\subset \M$ be such that $(f,\Lambda)$ is conjugated to a subshift of 
finite type. Then the result of Theorem~\ref{thm:1} holds.
\end{thm}

The proof of this result is decomposed in several parts.
First we give a construction similar to the one by Feng and Wu of points
with prescribed recurrence in a symbolic space \cite{fw}.

\subsection{Producing points with given symbolic recurrence sequence}
For a given sequence of integer $(\ell_n)$ such that

(a) $\exists n_0, \ell_{n+1}\ge \ell_n + 2n$ for any $n\ge n_0$ and

(b) $\lim_{n\to\infty} \ell_n/n^2=\infty$,

\noindent
we define a function $g$ nearly as in \cite{fw} such that points in the image
of $g$ have $(\ell_n)$ as recurrence sequence. Let $\A\subset\B$ be two alphabets
(i.e. finite or countable sets) with $\A\neq \B$ and $\sharp \B\ge 3$.
Choose a \emph{marker} $m\in\B\setminus \A$.

\begin{defi}\label{def:g}
Let $c\neq \bar c\in \B\setminus\{m\}$.
Define $g\colon \A^\NN \to \B^\NN$ by
$g(w)=\lim w^{(k)}$, where the sequence $w^{(k)}$ is constructed recursively as follows.

We put $w^{(0)}=w^{(1)}=\cdots=w^{(n_0-1)}=m w_1 w_2\cdots$.
For $k\ge n_0$ we construct $w^{(k)}$ from $w^{(k-1)}$ by inserting the
block $w^{(k-1)}_1\cdots w^{(k-1)}_k y_k$ at the position $\ell_k+1$, i.e.
\[
w^{(k)} = w^{(k-1)}_1 w^{(k-1)}_2 \cdots w^{(k-1)}_{\ell_k}
w^{(k-1)}_1\cdots w^{(k-1)}_k y_k w^{(k-1)}_{\ell_k+1} w^{(k-1)}_{\ell_k+2} \cdots,
\]
where $y_k=c$ if $w^{(k-1)}_{k+1}\neq c$ or $y_k=\bar c$ if $w^{(k-1)}_{k+1}=c$, so that
$y_k \neq m,w^{(k-1)}_{k+1}$.
Note that all the $w^{(k)}$ share the same $\ell_{n+1}$ first letters whenever $k\ge n$,
hence the limit $g(w)$ is well defined.
\qed
\end{defi}

Let us define the $k$-repetition time of $w$ by
\[
R_k(w) = \inf \{n>0\colon w_{n+1}w_{n+2}\cdots w_{n+k} = w_{1}w_{2}\cdots w_{k} \}.
\]
The main interest of the function $g$ resides in the following
\begin{lem}[\cite{fw}]\label{lem:g}
If $w\in \A^\NN$ and $k\ge n_0$ then $R_k(g(w))=\ell_k$.
\end{lem}
\begin{proof}[Proof of Lemma~\ref{lem:g}]
We do it by induction. Let $w^*=g(w)$.
Since there is no $m$ in $w$ it is obvious that $R_{n_0}(w^*)=\ell_{n_0}$.
Assume that for some $k\ge n_0$ we have proven that $R_k(w^*)=\ell_k$.
By construction this implies that $\ell_k < R_{k+1}(w^*)\le \ell_k+k$ or
$R_{k+1}(w^*)=\ell_{k+1}$
(note that there can be no $m$ between the positions $\ell_k+k+1$ and $\ell_{k+1}$ in $w^{*}$).
Suppose that we are in the first case and let $p=R_{k+1}(w^*)$.
The number of $m$ in the block $w^*_{p+1}\cdots w^*_{p+k} = w^*_1\cdots w^*_k$ being equal
to the number of $m$ inside the block $w^*_{\ell_k+1}\cdots w^*_{\ell_k+k}$ we get a contradiction
if $p>\ell_k$ (because $w^*_{\ell_k+1}=m$ is missed and there is no $m$ in
$w^*_{\ell_k+k+1}\cdots w^*_{p+k}$).
\end{proof}

Although one can imagine similar constructions for subshifts of finite type,
we decided to work with fullshifts to make the core argument more transparent.
However, we then have to work a bit to apply it to our systems, which are
presumably only coded by subshifts of finite type. This is the aim of next section.

\subsection{Building the source of large dimension}
Suppose that $\dim_H\Lambda>0$ otherwise there is nothing to prove
in Theorem~\ref{thm:2}. Let $(\Sigma,T)$ be a subshift of finite type
conjugated to $(f,\Lambda)$ and denote by $\pi\colon \Sigma\to\Lambda$ the conjugating homeomorphism.
For our purpose there is no loss of generality to suppose that $(\Sigma,T)$ is
topologically mixing.
Thus at least one letter, say $a$, has two successors
$b$ and $c$. Let then $B$ and $C$ be the shortest paths (and smallest in the lexicographic order)
starting in $b$ and $c$, respectively, and ending in $a$. That is $aB$ and $aC$ are the smallest
words of $\Sigma$ with prefix $ab$ and $ac$, respectively, and suffix $a$. It is clear that
there is a one-to-one correspondence between the set of words and the set of cylinders.
Without confusion, we also use a word to denote it's corresponding cylinder.
Let $A=aB$ and denote by $t(\omega)$ the return time of $\omega\in A$ into $A$,
\[
t(\omega)=\inf\{t>0\colon T^t\omega\in A\}.
\]
Note that $t$ is unbounded on $A$ since there is a disjoint union contained in $A$.
More precisely,
\[
A \supset aBB \cup aBCB \cup aBCCB \cup \cdots.
\]
To simplify the exposition we will consider the new alphabet consisting in
$|A|$-cylinders, where $|A|$ denotes the length of the cylinder $A$, together
with its associated transition matrix.
Let then $\Z_n$ be the partition of $\Sigma$ by $n$-cylinders (the new ones), and
denote by $\F_n$ the finite $\sigma$-algebra generated by $\Z_n$.

Let $\mu$ be the equilibrium state of the potential
$\varphi = -\dim_H(\Lambda) \log |Df|\circ\pi$. It is well known that
$\mu$ is the measure of maximal dimension and the Bowen pressure formula gives
$P(\varphi|\Sigma)=0$. In addition the measure $\mu$ is $\psi$-mixing (See e.g. \cite{k}
for properties of equilibrium states), which means that
\[
\psi(n) \eqdef\sup_{m} \sup_{U\in\F_m,V\in \B} \left\vert
\frac{\mu(U\cap T^{-n-m}V)}{\mu(U)\mu(V)}-1\right\vert \mathop{\longrightarrow}_{n\to\infty} 0.
\]
Hence there exists $m_0$ such that
\[
\mu(A \cap T^{-m_0}(\Sigma\setminus A) \cap T^{-2m_0}(\Sigma\setminus A) \cap \cdots \cap
T^{-m m_0}(\Sigma\setminus A)) \le  \delta^m,
\]
where $\delta\eqdef \mu(\Sigma\setminus A)(1+\psi(m_0)) <1$.
Set $H_n=\{ \omega\in A\colon t(\omega) \ge n\}$. The preceding estimate yields
\[
\mu(H_n) \le \delta^{n/m_0}.
\]
Moreover we have $H_n = A\setminus \{x\in A\colon t(x)\le n-1\}\in \F_n$. From
Proposition~\ref{pro:2} (see Section~\ref{sec:no-osc} below), setting
$\Sigma_n=(\Sigma\setminus H_n)^\infty$, we get that
\begin{equation}\label{pression}
\lim_{n\to\infty} P(\varphi | \Sigma_n ) = P(\varphi | \Sigma) = 0.
\end{equation}
Let $\nu_n$ be an ergodic equilibrium state of $\varphi$ on $\Sigma_n$.
If $n$ is sufficiently large then $\nu_n(A)>0$ (otherwise $\nu_n$ would be supported
on the forward invariant set $(\Sigma\setminus A)^\infty$, on which the pressure of $\varphi$ is strictly less than $P(\varphi | \Sigma)$).

Let $\hat \Sigma\subset A$ denotes the set of points returning infinitely many times in $A$.
Let us define the induced system $(\hat \Sigma, \hat T, \hat\nu_n)$ by
\[
\hat T(x)=T^{t(x)}(x), \quad\hat\nu_n = \frac {1}{\nu_n(A)}\nu_n|_A.
\]
Since $\hat T$ is the induced transformation on a cylinder of a subshift of finite type,
$(\hat \Sigma,\hat T)$ is a full shift with the countable alphabet
\[
\hat\Z = \bigcup_{m=1}^{\infty} \{ Z\cap\hat \Sigma\colon Z\in\Z_{m+1},Z\subset A,t|_Z=m \}.
\]
Since the return times in $A$ are bounded (by $n-1$) on $\Sigma\setminus H_n$, one has
\[
\hat \Sigma_n \eqdef A\cap \Sigma_n \subset \hat \Sigma.
\]
Observe that $(\hat \Sigma_n,\hat T)$ is still a fullshift over the finite alphabet
\[
\hat\Z^n = \bigcup_{m=1}^{n-1} \{ Z\cap\hat \Sigma_n\colon Z\in\Z_{m+1},Z\subset A,t|_Z=m \}.
\]
Set $\psi=\log|Df|\circ \pi$, $\displaystyle S_k\psi(\omega)=\sum_{j=0}^{k-1}\psi(T^j\omega)$
and $\displaystyle\hat S_k t(\omega)=\sum_{j=0}^{n-1}t(\hat T^j\omega)$.
Let us define \emph{the source of large dimension} by
\[
G_n= \left\{ \omega\in \hat \Sigma_{n} \colon
(i)\ \lim_{k\to\infty}\frac{1}{k}\ S_k\psi(\omega)\to \nu_n(\psi),
\, (ii)\ \lim_{k\to\infty}\frac{1}{k}\hat S_kt(\omega)\to \hat\nu_n(t)
\right\}.
\]
By Birkhoff's Ergodic Theorem we have $\hat\nu_n(G_n)=1$ which provides the following
lower bound
\begin{equation}\label{lwbd}
\dim_H \pi G_n\ge \dim_H \pi\hat\nu_n = \dim_H\pi\nu_n|_A = \dim_H\pi\nu_n.
\end{equation}
The last equality follows by ergodicity of the equilibrium state $\nu_n$.
Moreover,
\[
\dim_H\pi\nu_n = \frac{h_{\nu_n}}{\int\psi d\nu_n}
= \dim_H\Lambda + \frac{P(\varphi|\Sigma_n)}{\int\psi d\nu_n}.
\]
Observe that this together with \eqref{pression} and \eqref{lwbd} imply that
\begin{equation}\label{sourceld}
\lim_{n\to\infty} \dim_H \pi G_n = \dim_H \Lambda.
\end{equation}
\subsection{The image of the source has good symbolic recurrence}
With the alphabets $\A=\hat\Z^n$ and $\B=\hat\Z$ we can use Definition~\ref{def:g}
to produce a function $g\colon \hat \Sigma_n\to \hat \Sigma$.
For the marker $m\in\B\setminus\A$ we choose a $m\in\hat\Z$ such that $t|_m=n$.

Given an integer $k$, let $p$ be such that $\ell_p\le k<\ell_{p+1}$
and set $\epsilon_k = (p+2)^2/\ell_p$. The property (b) of the sequence $(\ell_p)$
implies that $\epsilon_k\to 0$ as $k\to\infty$.

Observe that $t$ is $\hat\Z$ measurable, hence for any $\omega\in G_n$ we have
\begin{equation}\label{eq:tg}
|\hat S_kt(g\omega) - \hat S_kt(\omega)| \le k\epsilon_k n,
\end{equation}
since the symbolic sequence of $g(\omega)$ and $\omega$ differs only because
of the block of length $i$ inserted at position $\ell_i$, which makes
a difference of at most $\sum_{i=n_0+1}^{p+1}i \le k\epsilon_k$ letters.
Hence for any $\omega\in G_n$ we have
\begin{equation}\label{averagetime}
\lim_{k\to\infty} \frac 1k \hat S_k t(g\omega) = \hat\nu_n(t).
\end{equation}

Notice that a $k$-cylinder of $(\hat \Sigma,\hat T)$ containing $\omega$ is indeed a
$\hat S_kt(\omega)$ cylinder of $(\Sigma,T)$.
Hence Lemma~\ref{lem:g} gives that for $\omega\in G_n$ we have
(the notation $\hat R_k$ stands for the $k$-repetition time computed with $\hat T$)
\[
R_{\hat S_kt(g\omega)}(g\omega)=\hat R_k(g\omega) = \ell_k.
\]
Passing to the limit using Equation~\eqref{averagetime} and the fact that the sequence $(R_h)_h$
is monotone gives us that the image $gG_n$ has good symbolic recurrence.
Namely, for any $\omega\in G_n$,
\begin{equation}\label{symbolicrec}
\liminf_{h\to\infty} \left(\limsup_{h\to\infty}\right) \frac {\log R_h(g\omega)}{h} =
\liminf_{k\to\infty} \left(\limsup_{k\to\infty}\right) \frac{1}{\hat\nu_n(t)}\frac {\log \ell_k}{k}.
\end{equation}

\subsection{The image has large dimension}
The fact that the dimension of the image $\pi gG_n$ is large will follow
from \eqref{sourceld} if we are able to prove that the inverse function of
$\gamma\eqdef\pi\circ g\circ \pi^{-1}$ has some smoothness.

Let $\omega\in G_n$. Denotes the $h$-cylinder of $(\Sigma,T,\Z)$
containing $g(\omega)$ by $\Z(g\omega,h)$.
If $k$ is such that
\[
\hat S_kt(g\omega)\le h < \hat S_{k+1}t(g\omega)
\]
then
\[
\Z(g\omega,h) \subset \hat \Z(g\omega,k),
\]
where $\hat\Z(z,k)$ denotes the $k$-cylinder of $(\hat \Sigma,\hat T,\hat \Z)$ containing a point $z$.
By construction of $g$ we have
\[
g^{-1}\hat\Z(g\omega,k) \subset \hat\Z(\omega,k(1-\epsilon_k)).
\]
Because of \eqref{eq:tg} we have
\[
\hat\Z(\omega,k(1-\epsilon_k)) \subset \Z(\omega,h(1-\delta_h))
\]
for some sequence $\delta_h\to 0$ as $h\to\infty$ (independent of $\omega$).
Putting together this chain of inclusions gives us that
\begin{equation}\label{g-cyl}
g^{-1} \Z(g\omega,h) \subset \Z(\omega,h(1-\delta_h)).
\end{equation}
This estimate on cylinders may be exploited for balls using the conjugacy assumption.
\begin{prop}\label{pro:5}
Let $(f,\Lambda)$ be a $C^{1+\alpha}$ conformal expanding map conjugated to a subshift.
Then there exists a constant
$\kappa\in(0,1)$ such that for any integer $n$ and any $x\in\Lambda$
\begin{equation}\label{ballcyl}
B(x, \kappa |D_x f^n|^{-1})\cap\Lambda \subset \pi\Z(\pi^{-1}x,n) 
\subset B(x,\kappa^{-1} |D_x f^n|^{-1})\cap\Lambda.
\end{equation}
\end{prop}
\begin{proof}
We give only a sketch. Assume that $\diam\Lambda=1$ and $|D_xf|>1$ for any $x\in\Lambda$.
Since $f$ is conjugated to a subshift the minimal distance between two image
$\pi Z$ and $\pi Z'$ of different elements $Z,Z'\in\Z$ is bounded from below
by some constant $\delta>0$.
Moreover, by distortion there exists a constant $D$ such that
for any $n$ and $x,y$ in the same $n$-cylinder we have
$|D_x f^n| \le D |D_y f^n|$.
Let $\kappa=\delta/D$.

Suppose that $d(x,y)<\kappa |D_xf^n|^{-1}$.
This implies that $d(x,y)<\delta$, thus $\Z(\pi^{-1}x,1)=\Z(\pi^{-1}y,1)$.
If $\Z(\pi^{-1}x,k-1)=\Z(\pi^{-1}y,k-1)$ for some $k<n$, then we get
$d(f^kx,f^ky) \le D |D_xf^k| d(x,y) < \delta$. Hence $\Z(\pi^{-1}x,k)=\Z(\pi^{-1}y,k)$.
By induction we get the first inclusion in Equation~\eqref{ballcyl}.
The remaining inclusion is easier and we omit its proof.
\end{proof}

\begin{lem}\label{lem:holder}
The function $\gamma^{-1}$ is $\alpha$-H\"older continuous on $g(G_n)$ for any $\alpha\in(0,1)$.
\end{lem}
\begin{proof}
By arguments similar to those in the proof of Equation~\eqref{eq:tg},
but using the continuity of $\psi$, we get that
\begin{equation}\label{lyap}
\lim_{n\to\infty} \sup_{\omega\in G_n} \
\left|\frac 1n S_n\psi(\omega) - \frac 1n S_n\psi(g\omega) \right| = 0.
\end{equation}
Let $x,y\in \pi G_n$. Let $n$ be such that
\[
\kappa \exp(-S_{n+1}\psi(g\pi^{-1}x))
\le d(\gamma x,\gamma y) < \kappa \exp(-S_{n}\psi(g\pi^{-1}x)).
\]
It follows from Proposition~\ref{pro:5} that $g\pi^{-1}y\in\Z(g\pi^{-1}x,n)$.
Thus $\pi^{-1}y\in \Z(\pi^{-1}x,n(1-\delta_n))$ by \eqref{g-cyl}.
Using Proposition~\ref{pro:5} again we find that
\[
d(x,y) < \kappa^{-1} \exp(-S_{n(1-\delta_n)}\psi(\pi^{-1}x)).
\]
For any $\alpha\in(0,1)$, because of \eqref{lyap}, the facts that $\delta_n\to0$,
$\psi$ is strictly positive and bounded, we can find an $n(\alpha)$ such that for all $n\ge n(\alpha)$ and all $x$
we have
\[
\exp(-S_{n(1-\delta_n)}\psi(\pi^{-1}x)) \le \kappa^2 \exp(-\alpha S_{n+1}\psi(g\pi^{-1}x))
\]
This shows that whenever $d(\gamma x,\gamma y)< \kappa\inf\exp S_{n(\alpha)}\psi$ we have
\[
d(x,y)< d(\gamma x,\gamma y)^\alpha,
\]
proving the lemma.
\end{proof}
Since $\pi G_n = \gamma( \pi g G_n)$ Lemma~\ref{lem:holder} ensures that the image of the
source has large dimension:
\begin{equation}\label{imageld}
\dim_H\pi g(G_n) \ge \dim_H \pi G_n.
\end{equation}

\subsection{The image has good recurrence rates}
By Proposition~\ref{pro:5} we have
\begin{equation}\label{tt}
\tau_{\kappa \exp(S_k \psi(\pi^{-1}x))}(x) \ge R_k(\pi^{-1}x) \ge \tau_{\kappa^{-1}\exp(S_k \psi(\pi^{-1}x))}(x),
\end{equation}
for any $x\in\Lambda$.
Let $\lambda_n = \int\psi d\nu_n$.
By Equation~\eqref{lyap}, $g$ does not change the Birkhoff average of $\psi$
thus for any $x\in \pi g(G_n)$, by the property (ii) in the definition of $G_n$
we have
\begin{equation}\label{lyapg}
\lim_{k\to\infty} \frac 1k S_k\psi(\pi^{-1}x) = \lambda.
\end{equation}
Then we choose the sequence $(\ell_k)$ as in \cite{fw} such that
\[
\liminf_{k\to\infty} \frac {\log \ell_k}{k} = \alpha\lambda_n\hat\nu_n(t)
\quad\textrm{and}\quad
\limsup_{k\to\infty} \frac {\log \ell_k}{k} = \beta\lambda_n\hat\nu_n(t).
\]
Using the bounds in \eqref{tt} and \eqref{lyapg} we get that for any $x\in \pi g(G_n)$
\[
\liminf_{r\to 0}     \frac{\log\tau_r(x)}{-\log r} =
\liminf_{k\to\infty} \frac{\log R_k(\pi^{-1}x)}{\lambda_n k} =
\liminf_{k\to\infty} \frac{1}{\lambda_n\hat\nu_n(t)}\frac{\log \ell_k}k =\alpha,
\]
where the second equality follows from \eqref{symbolicrec}.
The same arguments for the $\limsup$ gives that $x\in E(\alpha,\beta)$, where
\begin{equation}\label{eab}
E(\alpha,\beta) = \left\{x\in \Lambda\colon
\liminf_{r\to0}\frac{\log\tau_r(x)}{-\log r} = \alpha
\text{ and }
\limsup_{r\to0}\frac{\log\tau_r(x)}{-\log r} = \beta
\right\}.
\end{equation}
Since $\pi g(G_n)\subset E(\alpha,\beta)$, Theorem~\ref{thm:2} follows
from \eqref{sourceld} and \eqref{imageld}.

\section{Reduction to a system conjugated to a subshift}~\label{sec:no-osc}

We want to ensure that the original system contains subsystems with
arbitrary large dimension and which are conjugated to subshifts of finite type.
By Bowen's formula this is the same as finding such subsystems with
pressure of the potential $-\dim_H(\Lambda) \log|Df|$ arbitrarily close to $0$.
Fernandez, Ugalde and Urias prove that for maps of the interval with a Markov partition
by intervals this is always possible \cite{fuu}. Their proof relies on combinatorial
arguments. We provide here a proof based on the spectral
approach, close in spirit to the work by Maume and Liverani \cite{lm}. 
Although more intricate, this approach allows us to consider easily higher dimensional 
conformal maps.

\subsection{Pressure for systems with small hole}
Let $(X,T)$ be a subshift of finite type.
Given an open set $H\subset X$, the hole, we set $Y=X\setminus H$ and denote by
$Y^\infty=\cap_{m\ge 0} T^{-m}Y $ the compact forward invariant set of points never falling
into the hole $H$. Our next result says that the pressure of $Y^\infty$ tends to the
pressure of the original system as the hole gets smaller. Notice that this is not
a small perturbation for the potential, so that classical continuity of the pressure
do not apply.

\begin{prop}\label{pro:2}
Let $(X,T)$ be a topologically mixing subshift of finite type.
Assume that $\varphi\colon X\to\RR$ is a H\"older continuous potential.
Let $\mu$ be the equilibrium measure of the potential $\varphi$.
If $H_n$ is a decreasing sequence of $\F_n$ measurable sets such that
for some $c,\eps>0$
\[
\mu ( H_n ) \le c e^{-\eps n}
\]
then $\displaystyle\lim_{n\to\infty} P(\varphi |(X\setminus H_n)^\infty) = P(\varphi|X)$.
\end{prop}

\begin{proof}
The proof uses the characterization of the pressure by the logarithm of the spectral radius of
the Perron-Frobenius operator associated to the system. Then we show that one can
apply a perturbation result by Keller and Liverani \cite {kl} designed for these type of
operators.
Taking if necessary the new potential $\varphi-P(\varphi|X)+\psi\circ T-\psi$ we can assume
that $P(\varphi|X)=0$ and $L1=1$, where
$L$ is the Perron-Frobenius operator associated to $(X,T,\varphi)$, i.e.
\[
L f (x) = \sum_{y,Ty=x} \e^{\varphi(y)} f(y).
\]
Let $\Z_n$ be the partition of $X$ by $n$-cylinders, and denote by $\F_n$ the
finite $\sigma$-algebra generated by $\Z_n$.
By the Gibbs property of $\mu$ there exists $c_0$ such that for any $x\in Z\in\Z_k$
we have $1/c_0 \le \mu(Z)\e^{-\varphi_k(x)} < c_0$. We furthermore suppose
that for any cylinders $A,Z$ such that $AZ\eqdef A\cap T^{-k}Z\neq\emptyset$.
we have $1/c_0\le \mu(A)\mu(Z)/\mu(AZ) \le c_0$.

Given a function $f\colon X\to\RR$ and a set $A\subset X$ let $\osc(f,A)=\sup_A f -\inf_A f$.
Since $\varphi$ is H\"older there exists $c_1$ and $\alpha<1$ such that we have
\begin{equation}\label{eq:hold}
\var_n(\varphi)\eqdef \sup_{Z\in\Z_n}\osc(\varphi,Z) \le c_1 \alpha^n,\quad\text{for any $n\ge 1$.}
\end{equation}

Fix $\theta\in(1,\min(\alpha^{-1},\e^{\eps/2}))$. Define
\[
\osc(f)= \sum_{n\ge 1} \theta^n \sum_{Z\in\Z_n} \osc(f,Z)\mu(Z).
\]
Let $|f|_2=(\int |f|^2d\mu)^{1/2}$ and consider the Banach space
\[
B = \{ f\colon X\to\RR \colon \|f\|\eqdef\osc(f)+|f|_2 < \infty\}.
\]
We postpone the proofs of Lemmas~\ref{lem:2} and~\ref{lem:3} below to Section~\ref{sec:tec}.
\begin{lem}\label{lem:2}
There exists constants $c_3$ and $c_4$ such that for any $k\ge 1$
\[
\osc(L^kf) \le c_3\theta^{-k} \osc(f) + c_4 |f|_2.
\]
\end{lem}

Let $H_n^k = H_n \cup T^{-1}H_n\cup\cdots T^{-k+1}H_n$.
For simplicity we denote by $A$ instead of $1_A$ the indicator function of a set $A$.

\begin{lem}\label{lem:3}
There exists a constant $c_5$ (depending on $k$) such that for any $N\in\NN$ and $f\in B$,
\[
\osc(f H_N^k) \le \osc(f) + c_5 |f|_2.
\]
\end{lem}

For a set $Y\subset X$ we let $L_Y$ be the perturbed operator, $L_Yf=L(Yf)$.
Observe that $L_Y^kf = L^k(Y_kf)$, where $Y_k=Y \cap T^{-1}Y\cap\cdots T^{-k+1}Y$.
Hence $L_{X\setminus H_N}^kf= L^kf - L^k(H_N^k f)$.

\begin{lem}\label{lem:unicont}
There exists $k$ and $c_6$ such that for any $N\in\NN$ and $f\in B$ we have
\[
\osc( L_{X\setminus H_N}^kf) \le  \frac12 \osc(f) + c_6 |f|_2.
\]
\end{lem}
\begin{proof}[Proof of Lemma~\ref{lem:unicont}]
Choose $k$ such that $c_3\theta^{-k}\le\frac14$.
We apply Lemma~\ref{lem:2} and Lemma~\ref{lem:3} successively
\[
\begin{split}
\osc (L_{X\setminus H_N}^kf)
&\le \osc (L^kf)+ \osc(L^k(H_N^kf)) \\
&\le \frac14 \osc(f)+ c_4|f|_2 + \frac14\osc(H_N^kf) + c_4|H_N^kf|_2 \\
&\le \frac12 \osc(f)+ (2c_4+c_5) |f|_2.
\end{split}
\]
\end{proof}
The following lemma express that the hole is a perturbation of the original
system as far as operators from $B$ to $L^2(\mu)$ are concerned.
\begin{lem}\label{lem:perturb}
For any $N\in\NN$ and $f\in B$ with $\|f\|=1$ we have
\[
| (L-L_{X\setminus H_N})f |_2 \le c_0 \mu(H_N)^\frac14.
\]
\end{lem}
\begin{proof}[Proof of Lemma~\ref{lem:perturb}]
Using H\"older inequality then the invariance of $L$ we get
\[
\int |L(H_Nf)|^2 d\mu \le \sup_X L|f| \int_{H_N}|f| d\mu
\le \mu(H_N)^\frac12 |f|_2 \sup_X L|f|.
\]
Additionally, using again H\"older inequality, the supremum is bounded by
\[
\sup_X L|f| \le \sum_{Z\in\Z_1} c_0\mu(Z) \sup_Z |f|
\le \frac {c_0}\theta \osc(f) + c_0|f|_2.
\]
\end{proof}

It comes from the Folkolre Theorem that $ 1 = \e^{P(\varphi|X)}$ is a simple eigenvalue
of $L$ acting on $B$ and that the rest of the spectrum is contained in a disc of radius
$\rho<1$.
By Lemmas~\ref{lem:2},~\ref{lem:unicont}~and~\ref{lem:perturb} we can apply
Keller and Liverani's theorem (\cite{kl}, Corollary 1) with the sequence of
operators $L_{X\setminus H_n}$ acting on $B$. This gives that this spectral
figure will be conserved: for $n$ sufficiently large, the spectrum of
$L_{X\setminus H_n}$ on $B$ outside the disc of radius $\rho$ consists exactly
in one simple eigenvalue $\lambda_n$, with eigenvector $h_n$, such that
$\lambda_n\to 1$ and $\displaystyle h_n\cvto^B h$ (presumably $\lambda_n$ and $h_n$ are
positive, but we won't use this fact).

Furthermore, the restriction of the eigenvector $h_n$ to $(X\setminus H_n)^\infty$
is again an eigenvector of the Perron Frobenius operator $L_{(X\setminus H_n)^\infty}$
of the system $((X\setminus H_n)^\infty,T)$. Consequently the spectral radius of
$L_{(X\setminus H_n)^\infty}$ on $B$ is larger or equal to $|\lambda_n|$. Thus we get
$P(\varphi|X) \ge P(\varphi|(X\setminus H_n)^\infty) \ge \log |\lambda_n|$,
which proves the result.
\end{proof}

\subsection{Removing the boundary of the Markov partition}
The following proposition is of independent interest.
\begin{prop}\label{bdd}
Let $\Lambda\subset\M$ be a conformal repeller of the $C^{1+\alpha}$ map $f\colon \M\to \M $.
There exists subsets $\Lambda_n\subset \Lambda$ such that

(i) $(f,\Lambda_n)$ is conjugated to a subshift of finite type

(ii) $\displaystyle\lim_{n\to\infty} \dim_H\Lambda_n = \dim_H\Lambda$.
\end{prop}
\begin{proof}
Let $\xi$ be a generating Markov partition for the system $(f,\Lambda)$ and $(X,T)$
the corresponding subshift of finite type semi-conjugated via $\pi\colon X\to \Lambda$
to $(f,\Lambda)$. We write $\varphi = -\dim_H(\Lambda)\log|Df|$.

Let $K = \pi^{-1}\partial\xi$ be the premiage of the boundary of the partition.
The Markov property implies that $T K \subset K$. The compact invariant set $K$
is strictly contained in $X$ hence $P(\varphi | K)<0$ and the $\F_n$ measurable sets
$H_n = \Z_n(K)$ have a measure going exponentially to zero (see for example the
appendix of \cite{bs2} for a proof of these facts).
By Proposition~\ref{pro:2} for any $\epsilon>0$ there exists $n(\epsilon)$ such that
for any $n\ge n(\epsilon)$
\[
P(\varphi |(X\setminus H_n)^\infty ) > \epsilon\sup\varphi,
\]
whence the subset $\Lambda_n \eqdef \pi (X\setminus H_n)^\infty\subset \Lambda$ has dimension
\[
\dim_H(\Lambda_n) > (1-\epsilon) \dim_H(\Lambda).
\]
Moreover, $(f,\Lambda_n)$ is conjugated to a subshift of finite type.
\end{proof}
We are now able to give a proof of our main theorem.
\begin{proof}[Proof of Theorem~\ref{thm:1}]
By (i) in Proposition~\ref{bdd} we can apply Theorem~\ref{thm:2} to $(f,\Lambda_n)$.
This provides us with a set $E_n(\alpha,\beta)\subset\Lambda_n$ (See~\eqref{eab}) such that
$\dim E_n(\alpha,\beta)=\dim\Lambda_n$.
By monotony of Hausdorff dimension we get,
\[
\dim_H\Lambda \ge \dim_H E(\alpha,\beta) \ge \dim_H E_n(\alpha,\beta).
\]
The result follows from (ii) in Proposition~\ref{bdd}, taking the limit as $n\to\infty$.
\end{proof}

\section{Proofs of the technical lemmas}\label{sec:tec}

\begin{proof}[Proof of Lemma~\ref{lem:2}]
Let $\varphi_k = \varphi + \varphi\circ T + \cdots \varphi\circ T^{k-1}$.
Note that using Equation~\eqref{eq:hold} with $c_2=c_1/(1-\alpha)$, for any integer $n,k$, we have
$\var_{n+k}(\varphi_k) \le c_2\alpha^n$.
For a cylinder $A\in\Z_k$ and a point $x$ we denote by $Ax$ the unique element
of $A\cap T^{-k}x$, if any.
If $Z\in\Z_n$ with $n\ge 1$ then
\[
\begin{split}
\osc(L^kf,Z) &= \sup_{x,x'\in Z} L^kf(x)-L^kf(x')\\
&= \sup_{x,x'\in Z} \sum_{A\in \Z_k, AZ\neq\emptyset} e^{\varphi_k(Ax)} f(Ax) - e^{\varphi_k(Ax')} f(Ax')\\
&\le\sup_{x,x'\in Z}\sum_A e^{\varphi_k(Ax)} (f(Ax)-f(Ax'))
+ \e^{\varphi_k(Ax)}(\e^{\var_{n+k}(\varphi_k)}-1)\, |f|(Ax')\\
&\le \sum_A c_0\mu(A) \osc(f,AZ) + (\e^{\var_{n+k}(\varphi_k)}-1) c_0\mu(A)\sup_{AZ} |f|\\
&\le c_0(1+\e^{c_2}) \sum_A \mu(A) \osc(f,AZ) + \\
&\quad + c_0\e^{c_2}c_2\alpha^n \sum_A \mu(A)\mu(AZ)^{-1/2} \left(\int_{AZ}|f|^2d\mu\right)^{1/2}.
\end{split}
\]
Summing up over $Z\in\Z_n$ and then over $n$ yields to
\[
\begin{split}
\osc(L^kf)
&\le
\sum_{n\ge1}\theta^n\sum_{Z\in\Z_{n+k}} \Bigl[
c_0^2(1+\e^{c_2})\mu(Z) \osc(f,Z) + \\
& \quad + c_0^2\e^{c_2}c_2\alpha^n \mu(Z)^{1/2} \left(\int_{Z}|f|^2d\mu\right)^{1/2} \Bigr]\\
&\le c_3\theta^{-k} \osc(f) + c_4|f|_2 .
\end{split}
\]
For the last inequality we have set $c_3=c_0^2(1+\e^{c_2})$ and
$c_4=c_0^2\e^{c_2}c_2/(1-\alpha\theta)$ and use the Schwarz inequality.
\end{proof}

\begin{proof}[Proof of Lemma~\ref{lem:3}]
Since $H_N^k$ is $\F_{N+k}$ measurable we have
\[
\begin{split}
\osc(fH_N^k)
&= \sum_{n\ge 1} \theta^n \sum_{Z\in\Z_n} \osc(fH_N^k,Z)\mu(Z)\\
&\le \sum_{n=1}^{N+k-1} \theta^n \glu\sum_{Z\in\Z_n, Z\cap H_N^k\neq\emptyset}\glu
\mu(Z)\sup_Z|f|+ \sum_{n\ge N+k} \theta^n\sum_{Z\in\Z_n} \osc(f,Z)\mu(Z)\\
&\le \osc(f) + \sum_{n=1}^{N+k-1} \theta^n \glu\sum_{Z\in\Z_n, Z\cap H_N^k\neq\emptyset}\glu
\mu(Z)^{1/2} \left(\int_Z |f|^2d\mu\right)^{1/2} \\
&\le \osc(f) + \sum_{n=1}^{N+k-1} \theta^n\mu(\Z_n(H_N^k))^{1/2} |f|_2.
\end{split}
\]
For the last inequality we used Schwarz inequality and $\Z_n(H_N^k)$ denotes the
union of elements of $\Z_n$ intersecting $H_N^k$.
For any $n\le N+k-1$ we have
\[
\Z_n (H_N^k) \subset \cup_{j=0}^{k-1} T^{-j} \Z_{n-j}(H_N)
\subset\cup_{j=0}^{k-1} T^{-j} H_{n-k}.
\]
The result follows from the invariance of the measure by taking
\[
c_5 = \frac{\sqrt{ck\e^{\eps k}}}{1-\theta \e^{-\frac\eps2}}
\ge \sup_N \sum_{n=1}^{N+k-1} \theta^n (kc\e^{-\eps (n-k)})^{1/2} .
\]
\end{proof}

\subsection*{Acknowledgments}
The first author would like to thank Luis Barreira and Roman Hric for interesting
discussions on the topic and the organizers of the meeting ``Recent trends in Dynamics 2002'' held
in Porto, Portugal. The second author would like to thank CNRS UMR 6140, Amiens, where he
encountered warm hospitality.

\vfill
\end{document}